\documentclass[a4paper,12pt,reqno]{amsart}
\usepackage[a4paper,margin=1.9cm,top=2.5cm,bottom=2.5cm,centering,vcentering]{geometry}
\usepackage{amssymb,mathtools,cite,enumerate,color,eqnarray,hyperref}
\definecolor{blue}{RGB}{0, 0, 200}
\definecolor{pink}{RGB}{252, 0, 50}

\hypersetup{
           colorlinks=true,
           linkcolor=blue,
           urlcolor=red,
           citecolor=pink,
          }

\theoremstyle{plain}
\newtheorem{theorem}{Theorem}[section]

\newtheorem{lemma}[theorem]{Lemma}

\theoremstyle{definition}

\numberwithin{equation}{section}

\numberwithin{equation}{section}

\begin{document}

\title{Parity results of PEND partition}

\author[H. Nath]{Hemjyoti Nath}
\address{Department of mathematical sciences, Tezpur University, Napaam, Tezpur, Assam 784028, India}
\email{hemjyotinath40@gmail.com}


\date{\today.}

\begin{abstract}
In this paper, we consider the set of partitions $pend(n)$ which enumerates the number of partitions of $n$ wherein the even parts are not allowed to be distinct. Using a result of Newman, we prove a few infinite families of congruences modulo 2 for $pend(n)$.
\end{abstract}

\maketitle

\section{Introduction}
The partition of a positive integer $n$ is a non-increasing sequence of positive integers whose sum is equal to $n$. For example, $4+3+2+1$ is a partition of $10$.

If $p(n)$ denotes the number of partitions of a positive integer $n$, then with the convention that $p(0)=1$, the generating function of $p(n)$ (due to Euler) is given by
\begin{equation*}
    \sum_{n=0}^{\infty}p(n)q^n = \frac{1}{(q;q)_{\infty}},
\end{equation*}
where
\begin{equation*}
    (a;q)_{\infty} := \prod_{n=0}^{\infty}(1-aq^n),
\end{equation*}
where $a$ and $q$ are complex numbers with $|q|<1$. Throughout this paper, we set
\begin{equation*}
    f_k := (q^k;q^k)_{\infty}, \quad \text{for any integer} \quad k\geq1.
\end{equation*}

Ramanujan \cite{ram1}, \cite{ram2} and \cite{ram3} discovered three beautiful congruences for the partition function, namely
\begin{align*}
    p(5n+4)&\equiv 0 \pmod5, \\p(7n+5)&\equiv 0 \pmod7, \intertext{and} p(11n+7)&\equiv 0 \pmod{11}.
\end{align*}
The field grew rapidly as many mathematicians made important discoveries over the years. They found lots of new and interesting results. For example, they discovered many Ramanujan-Type congruences for different kinds of partitions. These include $l$-regular partitions, $t$-core partitions, and many more. Recently, Ballantine and Welch explored different ways to expand and refine POD and PED partitions, which led them to two new types of partitions they called POND and PEND partitions. In this paper, the focus is on PEND partitions, which are partitions where the even parts cannot be distinct.\\

The number of partitions of $n$ wherein the even parts are not allowed to be distinct is denoted by $pend(n)$. Sellers \cite{1}, proved that the generating function for $pend(n)$ is
\begin{equation}\label{e2}
     \sum_{n=0}^{\infty}pend(n)q^n = \frac{f_2f_{12}}{f_1f_4f_6}.
\end{equation}
He explored the $pend(n)$ function from an arithmetic perspective. In that same paper, he also discovered several congruences for $pend(n)$. For example, he proved that
\begin{equation*}
    pend(27n+19) \equiv 0 \pmod{3},
\end{equation*}
and
\begin{equation*}
    pend\left( 3^{2\alpha+1}n+\frac{17\cdot3^{2\alpha}-1}{8} \right) \pmod{3}, \quad \text{where} \quad \alpha \geq 1.
\end{equation*}

The main aim of this paper is to establish infinite families of congruences for the $pend(n)$ function. In the next theorem, we state our result.

\begin{theorem}\label{t3}
Let $a(n)$ be defined by
    \begin{equation*}
        \sum_{n=0}^{\infty}a(n)q^n = \frac{f_3^2}{f_1^3}.
    \end{equation*}
Let $p\geq 5$ be a prime and let $\left(\frac{\star}{p}\right)$ denote the Legendre symbol. Define
\begin{equation*}
    \omega(p) := a\left( \frac{p^2-1}{8} \right)+p^{-2}\left( \frac{-2}{p} \right)\left( \frac{\frac{-1}{8}(p^2-1)}{p} \right)
\end{equation*}
    \begin{enumerate}[(i)]
         \item If $pend\left( \frac{p^2-1}{8} \right) \equiv 1 \pmod{2}$, then for $n,k\geq 0$ and $1\leq j \leq p-1$, we have
    \begin{equation}\label{t3.1}
           pend\left( p^{4k+4}n + p^{4k+3}j + \frac{ p^{4k+4}-1}{8} \right) \equiv 0 \pmod{2},
        \end{equation}      
and for $k \geq 0$,
        \begin{equation}\label{t3.2}
            pend\left( \frac{p^{4k}-1}{8} \right) \equiv 1 \pmod{2}.
        \end{equation}

        \item If $pend\left( \frac{p^2-1}{8} \right) \equiv 0 \pmod{2}$, then for $n,k\geq 0$ and $1\leq j \leq p-1$, we have
      \begin{equation}\label{t2.2}
         pend\left( p^{6k+6}n + p^{6k+5}j + \frac{ p^{6k+6}-1}{8} \right) \equiv 0 \pmod{2},
         \end{equation}
and for $k\geq 0$, 
        \begin{equation}\label{t3.4}
            pend\left( \frac{p^{6k}-1}{8} \right) \equiv 1 \pmod{2}.
        \end{equation}

\end{enumerate}

\end{theorem}

\section{Preliminaries}
In this section, we recall some identities that will be used in our proofs. The well known Ramanujan's general theta function $f(a,b)$ \cite{6} is defined by
\begin{equation*}
    f(a,b)=\sum_{n=-\infty}^{\infty}a^{n(n+1)/2}b^{n(n-1)/2}.
\end{equation*}
Three special cases of $f(a,b)$ are the theta functions $\varphi(q)$, $\psi(q)$ and $f(-q)$, which are given by:
\begin{equation*}
    \varphi(q) := f(q,q) = \sum_{n=0}^{\infty}q^{n^2} = (-q;q^2)_{\infty}^2(q^2;q^2)_{\infty} = \frac{f_2^5}{f_1^2f_4^2},
\end{equation*}
\begin{equation*}
    \psi(q) := f(q,q^3) = \sum_{n=0}^{\infty}q^{n(n+1)/2} = \frac{(q^2;q^2)_{\infty}}{(q;q^2)_{\infty}} = \frac{f_2^2}{f_1},
\end{equation*}
and
\begin{equation*}
    f(-q):=f(-q,-q^2)=\sum_{n=0}^{\infty}(-1)^{n}q^{n(3n-1)/2} = (q;q)_{\infty}=f_1.
\end{equation*}
In terms of $f(a,b)$, Jacobi's triple product identity \cite{6} is given by
\begin{equation*}
    f(a,b) = (-a;ab)_{\infty}(-b;ab)_{\infty}(ab;ab)_{\infty}.
\end{equation*}

The following result of Newman will play a crucial role in the proof of our second theorem, therefore we shall quote it as a lemma. Following the notations of Newman's paper, we shall let $p$, $q$ denote distinct primes, let $r, s \neq 0, r \not \equiv s \pmod{2}$. Set
\begin{equation}\label{e7}
    \phi(\tau) = \prod_{n=1}^{\infty}(1-x^n)^r(1-x^{nq})^s = \sum_{n=0}^{\infty}c(n)x^n,
\end{equation}
$\epsilon = \frac{1}{2}(r+s), t=(r+sq)/24, \Delta= t(p^2-1), \theta = (-1)^{\frac{1}{2}-\epsilon}2q^s$, then the result is as follows:

\begin{lemma}[Newman \cite{16}]\label{l2}
    With the notations defined as above, the coefficients $c(n)$ of $\phi(\tau)$ satisfy
    \begin{equation}\label{e8}
        c(np^2+\Delta)-\gamma(n)c(n) +p^{2\epsilon-2}c\left( \frac{n-\Delta}{p^2} \right) = 0,
    \end{equation}
where 
\begin{equation}\label{e9}
    \gamma(n) = p^{2\epsilon - 2}\alpha-\left( \frac{\theta}{p} \right)p^{\epsilon-3/2}\left( \frac{n-\Delta}{p} \right),
\end{equation}
where $\alpha$ is a constant.
\end{lemma}

\section{Proof of theorem \ref{t3}}
We have
    \begin{equation}\label{e10}
        \sum_{n=0}^{\infty}a(n)q^n = \frac{f_3^2}{f_1^3}.
    \end{equation}
    Putting $r=-3, q=3$ and $s=2$ in $\eqref{e7}$, we have by Lemma $\ref{l2}$, for any $n\geq 0$
    \begin{equation}\label{e11}
        a\left( p^2n +\frac{p^2-1}{8} \right) = \gamma(n)a(n) - p^{-3}a\left( \frac{1}{p^2}\left( n-\frac{p^2-1}{8} \right) \right)=0,
    \end{equation}
    where
    \begin{equation}\label{e12}
        \gamma(n)=p^{-3}\alpha-p^{-2}\left(\frac{-2}{p}\right)\left(\frac{n-\frac{1}{8}(p^2-1)}{p}\right)
    \end{equation}
    and $\alpha$ is a constant integer. Setting $n=0$ in $\eqref{e11}$ and using the fact that $a(0) = 1$ and $a\left(\frac{\frac{-1}{8}(p^2-1)}{p^2}\right)=0$, we obtain
    \begin{equation}\label{e13}
       a\left( \frac{p^2-1}{8} \right) = \gamma(0).
    \end{equation}
Setting $n=0$ in $\eqref{e12}$ and using $\eqref{e13}$, we obtain
\begin{equation}\label{e14}
    p^{-3}\alpha = a\left( \frac{p^2-1}{8} \right)+p^{-2}\left(
    \frac{-2}{p}\right) \left(\frac{\frac{-1}{8}(p^2-1)}{p}\right) := \omega(p).
    \end{equation}
Now rewriting $\eqref{e11}$, by referring $\eqref{e12}$ and $\eqref{e14}$, we obtain
\begin{equation}\label{e15}
        a\left( p^2n +\frac{p^2-1}{8} \right) = \left( \omega(p)-p^2\left(\frac{-2}{p}\right)\left( \frac{n-\frac{p^2-1}{8}}{p} \right) \right)a(n) - p^{-3}a\left( \frac{1}{p^2}\left( n-\frac{p^2-1}{8} \right) \right).
    \end{equation}
Now, replacing $n$ by $pn + \frac{1}{8}(p^2-1)$ in $\eqref{e15}$, we obtain
\begin{equation}\label{e16}
    a\left( p^3n+\frac{p^4-1}{8} \right) = \omega(p)a\left( pn + \frac{p^2-1}{8} \right) - p^{-3} a(n/p).
\end{equation}

From equation $\eqref{e11}$, we can see that
\begin{equation}\label{e28}
        a\left( p^2n +\frac{p^2-1}{8} \right) - \gamma(n)a(n) + a\left( \frac{1}{p^2}\left( n-\frac{p^2-1}{8} \right) \right) \equiv 0 \pmod{2},
    \end{equation}
    where
    \begin{equation}\label{e29}
        \gamma(n)=\omega(p) + \left(\frac{n-\frac{p^2-1}{8}}{p}\right)
    \end{equation}
Setting $n=0$ in $\eqref{e28}$ and using the fact that $a(0)=1$ and $a\left( \frac{-\frac{p^2-1}{8}}{p^2} \right) = 0$, we arrive at
\begin{equation}\label{e30}
    a\left( \frac{p^2-1}{8} \right) = \gamma(0) \pmod{2}
\end{equation}
Setting $n=0$ in $\eqref{e29}$ yields
\begin{equation}\label{e31}
    \gamma(0) \equiv \omega(p) + 1 \pmod{2}.
\end{equation}
Combining $\eqref{e30}$ and $\eqref{e31}$ yields
\begin{equation}\label{e32}
    a\left( \frac{p^2-1}{8} \right) + 1 = \omega(p) \pmod{2},
\end{equation}

Again, from $\eqref{e2}$, we have that
\begin{equation}\label{e32.1}
    \sum_{n=0}^{\infty}pend(n)q^n = \frac{f_2f_{12}}{f_1f_4f_6} \equiv \frac{f_3^2}{f_1^3} = \sum_{n=0}^{\infty}a(n)q^n \pmod{2}
\end{equation}

\underline{\textbf{Case - 1 :} $a\left( \frac{p^2-1}{8} \right) \equiv 1 \pmod{2}$}

Since $a\left( \frac{p^2-1}{8} \right) \equiv 1 \pmod{2}$, so we have from $\eqref{e32}$ that, $\omega(p) \equiv 0 \pmod{2}$ and from equation $\eqref{e16}$ we obtain
\begin{equation}\label{e17}
    a\left( p^3n+\frac{p^4-1}{8} \right) \equiv p^{-3} a(n/p) \pmod{2}.
\end{equation}
Now, replacing $n$ by $pn$ in $\eqref{e17}$, we obtain
\begin{equation}\label{e18}
    a\left(p^4n+\frac{p^4-1}{8}\right) \equiv p^{-3} a(n) \equiv a(n) \pmod{2}.
\end{equation}
Since $p^{4k}n+\frac{p^{4k}-1}{8}=p^4\left(p^{4k-4}n+\frac{p^{{4k-4}}-1}{8}\right)+\frac{p^4-1}{8}$, using equation $\eqref{e17}$, we obtain that for every integer $k\geq 1$,
\begin{equation}\label{e19}
    a\left(p^{4k}n+\frac{p^{4k}-1}{8}\right) \equiv a\left(p^{4k-4}n+\frac{p^{4k-4}-1}{8}\right) \equiv a(n) \pmod{2}.
\end{equation}
Now if $p \nmid n$, then $\eqref{e17}$ yields
\begin{equation}\label{e20}
    a\left(p^3n +\frac{p^4-1}{8}\right) \equiv 0 \pmod{2}.
\end{equation}
Replacing $n$ by $p^3+\frac{p^4-1}{8}$ in $\eqref{e19}$ and using $\eqref{e20}$, we obtain
\begin{equation}\label{e21}
a\left( p^{4k+3}n + \frac{p^{4k+4}-1}{8} \right) \equiv 0 \pmod{2}.    
\end{equation}
In particular, for $1 \leq j \leq p-1$, we have from $\eqref{e21}$, that
\begin{equation}\label{e40}
    a\left( p^{4k+4}n + p^{4k+3}j + \frac{p^{4k+4}-1}{8}\right) \equiv 0 \pmod{2}
\end{equation}
Congruence $\eqref{t3.1}$ follows from $\eqref{e32.1}$ and $\eqref{e40}$. \\

Taking $n=0$ in $\eqref{e19}$ and utilizing the fact that $a(0)=1$, we get
\begin{equation}\label{e40.1}
    a\left(\frac{p^{4k}-1}{8}\right)  \equiv 1 \pmod{2}.
\end{equation}
Congruence $\eqref{t3.2}$ follows from $\eqref{e32.1}$ and $\eqref{e40.1}$. \\

\underline{\textbf{Case - 2 :} $a\left( \frac{p^2-1}{8} \right) \equiv 0 \pmod{2}$}

In order to prove $(ii)$, we replace $n$ by $p^2n+\frac{p(p^2-1)}{8}$ in $\eqref{e16}$
\begin{align}
    a\left( p^{5}n + \frac{p^6-1}{8} \right) & = a\left( p^3\left( p^2n+\frac{p(p^2-1)}{8} \right) + \frac{p^4-1}{8} \right) \nonumber\\
    & \equiv \omega(p)a\left(p^3n+\frac{p^4-1}{8}\right)-p^{-3}a\left( pn + \frac{p^2-1}{8} \right) \nonumber \\
    & \equiv \left[ \omega^2(p) - p^{-3}\right]a\left( pn+\frac{p^2-1}{8} \right) -p^{-3}\omega(p)a(n/p). \label{e22}
\end{align}
Now, as $a\left( \frac{p^2-1}{8} \right) \equiv 1 \pmod{2}$, so we have from  $\eqref{e32}$ that, $\omega(p) \equiv 1 \pmod{2}$ and $p\geq5$ is an odd prime, we have $\omega^2(p)-p^{-3}=0 \pmod{2}$, and therefore $\eqref{e22}$ becomes
\begin{equation}\label{e23}
    a\left( p^{5}n + \frac{p^6-1}{8} \right) \equiv a(n/p) \pmod{2}.
\end{equation}
Replacing $n$ by $pn$ in $\eqref{e23}$, we obtain
\begin{equation}\label{e24}
    a\left( p^6n +\frac{p^6-1}{8} \right) \equiv a(n) \pmod{2}.
\end{equation}
Using equation $\eqref{e24}$ repeatedly, we see that for every integers $k\geq 1$,
\begin{equation}\label{e25}
    a\left( p^{6k}n + \frac{p^{6k}-1}{8} \right) \equiv a(n) \pmod{2}.
\end{equation}
Observe that if $p\nmid n$, then $a(n/p)=0$. Thus $\eqref{e23}$ yields
\begin{equation}\label{e26}
    a\left( p^5n + \frac{p^6-1}{8} \right) \equiv 0 \pmod{2}.
\end{equation}
Replacing $n$ by $p^5n+\frac{p^6-1}{8}$ in $\eqref{e25}$ and using $\eqref{e26}$, we obtain
\begin{equation}\label{e27}
    a\left( p^{6k+5}n + \frac{p^{6k+6}-1}{8} \right) \equiv 0 \pmod{2}.
\end{equation}
In particular, for $1 \leq j \leq p-1$, we have from $\eqref{e27}$, that
\begin{equation}\label{e27.1}
    a\left( p^{6k+6}n + p^{6k+5}j + \frac{p^{6k+6}-1}{8}\right) \equiv 0 \pmod{2}
\end{equation}
Congruence $\eqref{t2.2}$ follows from $\eqref{e32.1}$ and $\eqref{e27.1}$. \\

Taking $n=0$ in $\eqref{e25}$ and utilizing the fact that $a(0)=1$, we get
\begin{equation}\label{e40.2}
    a\left(\frac{p^{6k}-1}{8}\right)  \equiv 1 \pmod{2}.
\end{equation}
Congruence $\eqref{t3.4}$ follows from $\eqref{e32.1}$ and $\eqref{e40.2}$. \\


\begin{thebibliography}{99}

\bibitem{2}
Ballantine C. and Welch A.,
PED and POD partitions : combinatorial proofs of recurrence relations.
\textit{Discrete Math.},
\textbf{346}, (2023).

\bibitem{3}
Ballantine C. and Welch A.,
Generalizations of POD and PED Partitions, preprint at https://arxiv.org/abs/2308.06136.


\bibitem{6}
Berndt, B.C., \textit{Number Theory in the Spirit of Ramanujan}, AMS, Providence, 2006.


\bibitem{5}
Hirschhorn, M.D., \textit{The power of \lowercase{$q$}, a personal journey}, 
Developments in Mathematics, 49.
Springer, 2017.

\bibitem{16}
Newman, M.
Modular forms whose coefficients possess multiplicative properties, $II$,
\textit{Ann. Math.},
\textbf{75}, (1962), 242--250.

\bibitem{ram1}
Ramanujan, S., Some properties of $p(n)$, the number of partitions of $n$, \textit{Proceedings of the Cambridge Philosophical Society}, 19:207–210, (1919).

\bibitem{ram2}
Ramanujan, S., Congruence properties of partitions, \textit{Proceedings of the London Mathematical Society}, 18:19, (1920).

\bibitem{ram3}
Ramanujan, S., Congruence properties of partitions, \textit{Mathematische Zeitschrift}, 9(1-2):147–153, (1921).

\bibitem{1}
Sellers, J.A.,
Elementary proofs of congruences for pond and pend partitions,
\textit{J. Integer Seq.},
\textbf{24} (2024).

\end{thebibliography}
\end{document}